# A Geometric Solution to the Isoperimetric Problem and its Quantitative Inequalities

By

Lakshya Chaudhary

Affiliation - Vidya Academy School, Daurala, (Sardhana)-250221, Meerut, Uttar Pradesh, India. (High School Student)

E-mail – lakshyachaudhary222@gmail.com

# 1. Abstract

This paper presents a geometric approach to the classical isoperimetric problem **[1]** by analysing the efficiency of regular polygons in enclosing maximum area for a fixed perimeter. Using efficiency metrics, it proves that regular polygons converge to the circle in area efficiency as the number of sides increases. The paper also extends these ideas to irregular polygons using average apothem and circumradius, defines "wasted area," and introduces several compactness and smoothness indices, offering a unified and quantitative view of isoperimetric inequalities **[2]**

# 2. Introduction

The isoperimetric problem is one of the oldest and most elegant questions in geometry: among all shapes with a given fixed perimeter, which one encloses the greatest possible area? The classical result is well known—the circle is the unique maximised area for a given perimeter. While this has was proven analytically using calculus of variations, and more recently refined through quantitative isoperimetric inequalities, intuitive geometric insights into this phenomenon remain as a rich area of study.

In this paper, we approach the isoperimetric problem from a constructive and geometric perspective by examining the efficiency of regular polygons as they approximate the circle. We introduce a new formula for polygonal efficiency based on the ratio of a polygon's apothem to the radial hypotenuse (the distance from the centre to a vertex), normalized by the perimeter. This ratio measures how effectively a given shape uses its perimeter to enclose area.

We derive a new expression for the area of a regular polygon in terms of its perimeter and apothem, and define a wasted area function representing the gap between the area of a polygon and the area of a circle with the same perimeter. We derived many efficiency metrics for regular polygons.

Furthermore, this idea can be extended to non-regular convex polygons by approximating the area using the average apothem and average circumradius. This provides a useful geometric tool for estimating the area of irregular shapes without relying on coordinate-based methods

## 2. A Geometric Solution to the Isoperimetric problem by Regular Polygonal Convergence and Wasted Area

The isoperimetric problem seeks to determine which plane figure of a given perimeter encloses the maximum area. The classical solution is the circle. In this section, we present a constructive geometric proof based on comparing the area of regular polygons with increasing number of sides to that of a circle, using the concept of area wastage.

**3.1. Regular Polygonal Setup –** Let $n \in \mathbb{N}, n > 2$ represent the number of sides in a regular polygon. Let be the fixed total perimeter 'p' of the regular polygon.

Each side has length : $\frac{p}{n}$

**3.2. Convergence of n-sided regular polygon to a circle –** For example, we have a n-sided regular polygon with perimeter 'p', if we increase its sides with constant perimeter, the shape will become more and more smooth with increasing sides and converge to a circle with, n = ∞ . With increasing sides, the area enclosed by the regular polygon will also increase because when the number of sides increases, the distance between the centre of the regular polygon and points on the boundary becomes less varied, as we choose different points on the boundary of that regular polygon. This led to an increase in the average distance between the points on the boundary and the centre of regular polygon. So, the coverage of area by the regular polygon also increases in all directions as it becomes perfectly smooth and with $n = \infty$, it becomes smoothest and occupies largest area with constant perimeter.

**3.3. Apothem and the Area of a Regular Polygon –** To calculate area of a regular polygon with 'n' sides in form of its perimeter, we need to divide the polygon into equal isosceles triangles. The number of triangles is equal to number of sides of that regular polygon. The triangle has one vertices touching the centre of the polygon and other two vertices are two ends of a side of the polygon. The base of triangle is the length of a side of the polygon and the height is the apothem (The distance from the centre of the polygon to the midpoint of a side.

Let $a_n$ *denote the apothem,*

$$a_n = \frac{s}{2\tan\left(\frac{\pi}{n}\right)}$$

$$= \frac{p}{2n \cdot \tan\left(\frac{\pi}{n}\right)}$$

The area of one triangle is :

$$A_{triangle} = \frac{1}{2} \cdot s \cdot a_n$$

The total area of the polygon is :

$$A_n = n \cdot A_{triangle} = \frac{1}{2} \cdot a_n \cdot p$$

$$= \frac{p^2}{4n \cdot \tan\left(\frac{\pi}{n}\right)}$$

While the area of a regular polygon is commonly derived using apothem-based methods, an alternative formulation uses the central angle and the circumradius of the polygon. This approach provides another geometric insight into how regular polygons approximate a circle as the number of sides increases.

Let $n > 2$, be the number of sides in a regular polygon with fixed perimeter 'p'. Each side of length $s = \frac{p}{n}$ subtends a centre angle :

$$\theta = \frac{2\pi}{n}$$

Let 'R' be the circumradius, i.e., the distance from the centre of the polygon to any vertex. Each of the isosceles triangles 'n' that make up the polygon has two sides of length and an included angle .

The area of each triangle can be computed using the trigonometric formula :

$$A_{triangle} = \frac{1}{2} R^2 \sin(\theta) = \frac{1}{2} R^2 \sin\left(\frac{\pi}{n}\right) \ \textbf{[3]}$$

So, the total area of the polygon is :

$$A_n = n \cdot A_{triangle} = \frac{n}{2} \cdot R^2 sin\left(\frac{2\pi}{n}\right)$$

To express this area in term of perimeter 'p', we need to use law of sine in one of the triangles :

$$\mathrm{s} = 2\mathrm{R} \cdot \sin\left(\frac{\pi}{\mathrm{n}}\right) \ \textbf{[4]} \Rightarrow \ \mathrm{R} = \frac{\mathrm{s}}{2 \sin\left(\frac{\pi}{\mathrm{n}}\right)} = \frac{\mathrm{p}}{2\mathrm{n} \times \sin\left(\frac{\pi}{\mathrm{n}}\right)}$$

Substitute this into the area formula :

$$A_n = \frac{n}{2}\left[\frac{p}{2n \cdot sin\left(\frac{\pi}{n}\right)}\right]^2 \cdot \sin\left(\frac{2\pi}{n}\right) = \frac{p^2}{8n \cdot \sin^2\left(\frac{\pi}{n}\right)} \cdot \sin\left(\frac{2\pi}{n}\right)$$

This expression is equivalent to the apothem-based area formula :

$$A_n = \frac{p^2}{4n \cdot \tan\left(\frac{\pi}{n}\right)}$$

This reinforces the geometric intuition by focusing on angular structure, highlights the symmetry and circular arc approximation and shows that multiple geometric pathways lead to same isoperimetric limit. It also paves the way for exploring non-Euclidean analogues or applying similar techniques to polyhedra in 3D, where circumradius-based definitions are more natural than apothem-based ones.

**3.4. Wasted Area Definition –** Suppose a regular polygon inscribed inside a circle with same perimeter. The regular polygon with finite sides always has less area than the circle which has same perimeter. The area between the boundary of the circle and polygon can be termed as 'wasted area', the more the sides of the regular polygon, the less the 'wasted area'. For

Let $a_{cirle}$ denote the area of a circle with the same perimeter 'p'. For a circle :

$$Perimeter = \ 2\pi r$$

Radius in terms of perimeter, $r = \frac{p}{2\pi}$

Thus, Area = $\pi r^2 = \pi\left(\frac{p}{2\pi}\right)^2 = \frac{p^2}{4\pi}$

Now, we have area of a circle and n-sided regular polygon in terms of their perimeter with same perimeter. So, we can find the 'wasted area'.

$$W_n = \ A_{circle} - \ A_n$$

$$= \frac{p^2}{4\pi} - \frac{p^2}{4n \cdot \tan\left(\frac{\pi}{n}\right)}$$

$$= \frac{p^2}{4}\left[\frac{1}{\pi} - \frac{1}{n \cdot \tan\left(\frac{\pi}{n}\right)}\right]$$

This represents how much area is "lost" when using a regular polygon instead of a circle.

**3.5. The Shape with the Largest Area Encloses with the Same Perimeter :** By 'section 3.2', we inferred that the circle is the shape that encloses the largest area with same perimeter i.e. The Shape with 100% area efficiency. To confirm this, we need to find the 'wasted area' with a infinite sided regular polygon.

To understand how $W_n$ behaves as $n \to \infty$, we expand the tangent function using the Taylor series **[5]**

$$A_n = \frac{p^2}{4n\left(\frac{\pi}{n} + \frac{\pi^3}{3n^3} + \cdots\right)}$$

$$\approx \frac{p^2}{4\pi}\left(1 - \frac{\pi^2}{3n^2} + \cdots\right)$$

Then $W_n$ will be,

$$W_n = \frac{p^2}{4\pi} - A_n \approx \frac{p^2}{4\pi} \cdot \frac{\pi^2}{3n^2}$$

$$= \frac{p^2\pi}{12n^2}$$

To find the 'wasted area' for infinite-sided regular polygon (circle), we need to set n = ∞.

So, 'wasted area' will be,

$$W_\infty = limit_{n\to\infty}(W_n)$$

$$= \ limit_{n\to\infty}\left[\frac{p^2\pi}{12n^2}\right] = 0$$

And,

$$A_\infty = \left(W_\infty + \frac{p^2}{4\pi}\right)\ or\ limit_{n\to\infty}(A_n)$$

$$= \frac{p^2}{4\pi}$$

We can also prove this using only basic trigonometric inequalities, that the area of a regular polygon is always less than the area of a circle with the same perimeter.

For a regular n-gon :

$$A_n = \frac{p^2}{4n \cdot \tan\left(\frac{\pi}{n}\right)} < \frac{p^2}{4\pi} = \ A_{circle}$$

Dividing both sides by $\frac{p^2}{4}$, this reduce to

$$n \cdot \tan\left(\frac{\pi}{n}\right) > \pi$$

This holds for all $n > 2,$ since $\tan(x) > x$ for $x \in \left(0, \frac{\pi}{2}\right)$. This also confirms the isoperimetric inequality through elementary means.

**3.6. Upper Bound on the Wasted Area –** While the exact expression for the area wasted by a regular polygon with perimeter  compared to a circle of the same perimeter has been derived earlier, it is often useful to establish a strict and simplified upper bound. This allows us to understand how the polygon approaches the circle quantitatively as the number of sides increases. We aim to bound $W_n$ from previously used trigonometric estimates.

We used the known inequality valid for all $x \in \left(1, \frac{\pi}{2}\right)$ :

$$\tan(x) > x + \frac{x^3}{3}$$

Setting $x = \frac{\pi}{n}$, we get

$$\tan\left(\frac{\pi}{n}\right) > \frac{\pi}{n} + \frac{\pi^3}{3n^3}$$

Taking reciprocals (preserving inequality direction since all terms are positive)

$$\frac{1}{\tan\left(\frac{\pi}{n}\right)} > \frac{1}{\frac{\pi}{n} + \frac{\pi^3}{3n^3}} = \frac{n}{\pi + \frac{\pi^3}{3n^2}}$$

Now Substitute into the area formula for the polygon for bounding the polygon area :

$$A_n = \frac{p^2}{4n \cdot \tan\left(\frac{\pi}{n}\right)} < \frac{p^2}{4} \cdot \frac{1}{\pi + \frac{\pi^3}{3n^2}}$$

Subtract this from circle's area :

$$W_n = A_{circle} - A_n < \frac{p^2}{4\pi} - \frac{p^2}{4\left(\pi + \frac{\pi^3}{3n^2}\right)}$$

Combine the fractions :

$$W_n = \frac{\pi p^2}{4(3n^2 + \pi^2)}$$

Thus, for all n > 2, the wasted area satisfies :

$$W_n = \frac{\pi p^2}{4(3n^2 + \pi^2)}$$

This is a valid and strict upper bound, which decreases with 'n' , and shows that the area lost due to angularity shrinks at least quadratically. This complements our earlier result :

$$W_n \approx \frac{p^2\pi}{12n^2}$$

and reinforces the convergence of the polygonal shape to the circle from a bounding perspective.

This upper bound confirms that the circle not only encloses the maximal area for a fixed perimeter, but that no regular polygon can surpass it. Moreover, the convergence of a polygon to the circular optimum is quantifiable and sharp, making this bounding method a useful analytical tool.

**3.7. Monotonic Increase of Polygon Area with Increasing Sides –** A natural question arising in the study of regular polygons and their convergence to the circle is whether the enclosed area strictly increases as the number of sides increases, for a fixed perimeter. In this subsection, we prove that the area of a regular n-gon strictly increases with 'n', supporting the geometric intuition that more sides result in a better approximation to a circle.

We know that the area of a regular polygon is :

$$A_n = \frac{p^2}{4n \cdot \tan\left(\frac{\pi}{n}\right)}$$

We need to prove that,

$$A_{n+1} > A_n \text{ for all } n > 2$$

To prove This, we can compare the area formula for successive values on 'n'. The inequality :

$$A_{n+1} > A_n$$

$$\frac{1}{\tan(n+1) \cdot \left(\frac{\pi}{n+1}\right)} > \frac{1}{n \cdot \tan\left(\frac{\pi}{n}\right)}$$

Multiplying both sides by $\tan\left(\frac{\pi}{n+1}\right) \tan\left(\frac{\pi}{n}\right). \mathit{We\ obtain}$ :

$$n \cdot \tan\left(\frac{\pi}{n}\right) > (n+1) \tan\left(\frac{\pi}{n+1}\right)$$

The inequality hold because the function $f(n) = n \cdot \tan\left(\frac{\pi}{n}\right)$ is strictly decreasing for $n > 2$. This follows from the fact that $\frac{\tan(x)}{x}$ is strictly increases for $x \in \left(0, \frac{\pi}{2}\right)$, which implies :

$$\frac{\tan\left(\frac{\pi}{n}\right)}{\frac{\pi}{n}} > \tan\left(\frac{\pi}{n+1}\right) \Rightarrow \; n \cdot \tan\left(\frac{\pi}{n}\right) > (n+1)\tan\left(\frac{\pi}{n+1}\right)$$

Hence, the area is strictly increasing with 'n'.

This result reinforces the idea that as the number of sides increases, a regular polygon becomes a more efficient shape for enclosing area under a fixed perimeter. Combined with earlier results showing convergence to the circle, this monotonicity confirms that the circle is not just optimal, but also the limit of a sequence of progressively improving polygonal shapes.

## 4. Polygonal Efficiency Metrics and Convergence

In this section, we introduce and analyse geometric efficiency metric that quantify how closely a regular polygon approximates the optimal shape—the circle—or a given perimeter. Unlike Section 2, which focused on area convergence, here we develop structural and normalized comparisons that offer additional insights.

**4.1. Smoothness of a polygon :** For a regular polygon with side length $s$ , apothem $a_n$ , and circumradius $R_n$, each triangle formed from the centre defines a smoothness measure : $Smoothness_n(\sigma_n) = \frac{a_n}{R_n}$

We interpret this as a smoothness ratio — it reflects how "round" or smooth the polygon is. As $n \to \infty, A_n \to R_n$, and this ratio tends to 1, matching the geometry of a circle. For smaller $n$, the polygon has sharp corners, and the apothem is significantly shorter than the radius, yielding a smaller ratio.

This ratio serves as a local geometric indicator: the closer it is to 1, the more structurally similar the polygon is to a circle in terms of smooth curvature.

We have :

$$A_n = \frac{p}{2n \cdot \tan\left(\frac{\pi}{n}\right)}$$

$$R_n = \frac{p}{2n \cdot \tan\left(\frac{\pi}{n}\right)}$$

Then,

$$Smoothness_n(\sigma_n) = \frac{a_n}{R_n} = \frac{\left[\frac{p}{2n \cdot \tan\left(\frac{\pi}{n}\right)}\right]}{\left[\frac{p}{2n \cdot \sin\left(\frac{\pi}{n}\right)}\right]} = \frac{\tan\left(\frac{\pi}{n}\right)}{\sin\left(\frac{\pi}{n}\right)}$$

$$= \cos\left(\frac{\pi}{n}\right)$$

smoothness % will be :

$$\sigma_n\% = \sigma_n \cdot 100$$

$$= 100 \cdot \cos\left(\frac{\pi}{n}\right)$$

**4.2. Perimeter Slope Efficiency :** In traditional isoperimetric analysis, perimeter is treated as fully horizontal for estimating enclosed area. However, in polygons—especially when viewed as composed of radial triangles—the boundary includes slanted sides, distributing perimeter vertically as well as horizontally. This area-enclosing efficiency. To quantify this, we define the Perimeter Slope Efficiency, which adjusts for the average slope '$m$' of the boundary triangles. The total vertical rise across the polygon is approximated as '$2am$' (where '$a$' is the apothem). Using the Pythagorean theorem, we combine this vertical component with the perimeter length to calculate the effective slope-corrected perimeter efficiency.

$$Effective\ Perimeter\ Length = \sqrt{p^2 + (2am)^2}$$

Then, the perimeter slope efficiency is defined as :

$$E_s = \frac{p}{\sqrt{p^2 + (2am)^2}}$$

If $m$ =0 (perfectly flat sides), Then :

$$E_s = \frac{p}{\sqrt{p^2}} = 1$$

If $m \to \infty$, the vertical dominates :

$$E_s \to \infty$$

The formula is dimensionless and bounded between 0 and 1. This makes the formula especially suitable for modeling how slope weakens a polygon's ability to approximate a circle — which is the ideal in isoperimetric geometry.

**4.3. Slant-Angle Efficiency Index :** To quantify the geometric efficiency of regular polygons with respect to their internal slant, we construct a compactness index involving the apothem, side slope, and perimeter. The aim is to produce a smooth, bounded function that shows how close a polygon is to circular behaviour based on the triangular slant geometry. We begin with the idea that a triangle inside a regular polygon, drawn from the centre to a side, stretches based on both the apothem and the slope of the side. A higher value of '$am$' indicates a steeper or more "slanted" structure. To measure this distortion relative to the overall size of the polygon, we normalize it using the perimeter '$p$', forming the dimensionless distortion ratio :

$$D = \frac{2am}{p}$$

Here, the constant '2' arises from considering the full triangle width from the centre to both endpoints of a side (symmetrically placed). Rather than using this distortion directly, we seek a bounded index that decreases smoothly with increasing distortion and lies strictly between 0 and A classical method to do this is applying the normalization function :

$$E = \frac{1}{\sqrt{1+D^2}}[\mathbf{6}] = \frac{1}{\sqrt{1+\left(\frac{2am}{p}\right)^2}}$$

As the product '$am$' decreases (i.e., the triangle becomes flatter and more circular), the efficiency index approaches 1, indicating high compactness. As '$am$' increases (sharper slant), efficiency decreases toward 0, reflecting lower geometric efficiency. The function is smooth, continuous, and strictly bounded in (0,1).

**4.4. Isoperimetric Slope Index :** Regular polygons inscribed in a circle are often analysed through global isoperimetric ratios such as $\frac{A}{p^2}$. However, such metrics do not fully account for local triangular geometry—specifically, how efficiently each constituent triangle within the polygon contributes to the total enclosed area based on its slant structure. To incorporate this internal geometry, we define the Isoperimetric Slope Index as a triangle-based area-to-slant ratio that generalizes compactness. Let a regular n-gon be composed of '$n$' congruent isosceles triangles, each with :

$$Base = s = \frac{p}{n} \; and \; half \; base = \frac{s}{2} = \frac{p}{2n}$$

$$Height(apothem) = \frac{p}{2n \cdot \tan\left(\frac{\pi}{n}\right)}$$

By the Pythagorean theorem, the slant length (hypotenuse) of each triangle is :

$$h_n = \sqrt{\left(\frac{p}{2n}\right)^2 + \left[\frac{p}{2n \cdot \tan\left(\frac{\pi}{n}\right)}\right]^2}$$

The area of the entire polygon, expressed in terms of '$p$' and '$n$', is :

$$E_{iso-slope} = \frac{\left[\frac{p^2}{4n \cdot \tan\left(\frac{\pi}{n}\right)}\right]}{\left(\frac{p}{2n}\right)^2 + \left[\frac{p}{2n \cdot \tan\left(\frac{\pi}{n}\right)}\right]^2}$$

We square the slant length  in the denominator to ensure dimensional consistency and shape-scale invariance. This mirrors classical compactness ratios such as , and guarantees that the Isoperimetric Slope Index remains a true geometric efficiency, independent of polygon size. This area expression derives from summing the areas of  congruent isosceles triangles. A polygon with short, wide triangles (small ) will yield a higher index, while tall, narrow triangles (large ) will reduce the index. As the polygon becomes circular, triangles flatten, and the slant length becomes minimal relative to area, maximizing the index. This makes this bounded Isoperimetric Slope Index a useful tool for comparing how effectively different polygons approach circular compactness from a triangle-centric perspective.

**4.5. Apothem-Angle Index :** This index measures how effectively a regular polygon's apothem reaches inward relative to how open each triangle is. As we divided the polygon into '$n$' congruent isosceles triangles, each with central angle $\theta = \frac{2\pi}{n}$. If we split one of these triangles in half, we obtain a right triangle with apothem '$a$' as one leg, hypotenuse '$h$' from the centre to a vertex, and angle $\frac{\pi}{n}$ at the centre. Using trigonometry, we have $\cos\left(\frac{\pi}{n}\right) = \frac{a}{h}$ which implies $h = \frac{a}{\cos\left(\frac{\pi}{n}\right)}$ at the centre. Using this, the difference between hypotenuse and the apothem gives a sense of how much geometric distance is "wasted" due to the triangle's angular spread. This inefficiency is modeled by the difference $1 - \cos\left(\frac{\pi}{n}\right)$ which approaches zero as the polygon becomes more circular. We multiply '$a$' by '$n$', and divide both by '$p$' and angular efficiency term $1 - \cos\left(\frac{\pi}{n}\right)$. This yields this unbounded expression :

$$E_{raw} = \frac{a \cdot n}{p\left[1 - \cos\left(\frac{\pi}{n}\right)\right]}$$

To normalize the index within (0,1), we apply the transformation :

$$E_{angle-apothem} = \frac{1}{1 + \ p\left[1 - \frac{\cos\left(\frac{\pi}{n}\right)}{a \cdot n}\right]}$$

This dimensionless and scale-invariant increase monotonically, as $n \rightarrow \infty$, this approaches 1.

**4.6. Radial Packing Efficiency :** A regular polygon inscribed in a circle approximates the circular boundary using straight-line segments. However, this leads to a measurable loss in radial coverage compared to a smooth circle. To quantify this, we define the Radial Packing Efficiency as the ratio of the polygon's area to the area of the circle in which it is inscribed.

$$A_c = \pi r^2$$

The area of each triangle is given by :

$$A_t = \frac{1}{2} R^2 \sin\left(\frac{2\pi}{n}\right)$$

Multiplying by the number of triangles gives the total polygon area :

$$E_r(n) = \frac{A_p}{A_c} = \frac{\frac{n}{2} R^2 \sin\left(\frac{2\pi}{n}\right)}{\pi r^2}$$

We can also derive this by using the total interior angle sum of a regular polygon, which is $(n-2)\pi$. Defining an angular efficiency based on how this angle sum contributes to area enclosure, we compare the polygon's area to $(n-2)\pi r^2$, yielding :

$$E_{angle}(n) = \frac{A_p}{(n-2)\pi r^2} = \frac{\mathrm{n} \cdot \sin\left(\frac{2\pi}{n}\right)}{2\pi(n-2)}$$

$$E_{angle}(bounded) = \frac{n \cdot \sin\left(\frac{2\pi}{n}\right)}{2\pi}$$

Therefore, both the radial fill argument and the interior angle normalization lead to the same compactness efficiency formula. This expression is dimensionless and lies strictly between 0 and 1. As the number of sides increases, the central angle becomes smaller and the polygon better approximates the circle. In the limit as $n \to \infty$, we have :

$$\lim_{n\to\infty} E_r(n) = 1$$

This reflects that the polygon fills the circle perfectly when it becomes infinitely sided. The formula thus provides a smooth and accurate measure of radial efficiency, complementing perimeter- and slope-based indices.

**4.7. Area Efficiency by Angle Sum and Sector Fill Ratio :** Regular polygons get closer and closer to circles as the number of their sides increases. But one way to measure how efficiently a polygon encloses area is to ask: how much area does it produce for the total amount of angular "turn" it contains? We can think of this angle sum as a kind of internal angular budget. Since area is measured in square units and angles in radians are dimensionless, we scale this angular budget by the square of the radius of the circle in which the polygon is inscribed. That gives us an angular reference area :

$$A_{ref} = (n-2)\pi R^2$$

This acts as a benchmark: if the polygon made perfect use of its angles, how much area would it be enclosing? We then compare the polygon's actual area to this reference. The result is the angle based area efficiency :

$$E_{angle} = \frac{A_n}{(n-2)\pi R^2}$$

This number tells us how effectively the polygon uses its angular structure to enclose space. It's dimensionless and depends only on the shape, not its size. As increases, the total angle sum grows, but the area approaches a fixed limit—namely, the area of the circle $\pi R^2$. That means this efficiency decreases with '$n$', and in the limit becomes :

$$\lim_{n\to\infty} E_{angle} = \frac{\pi R^2}{(n-2)\pi R^2} = \frac{1}{(n-2)} = 1$$

We could try to normalize this formula to make it bounded by multiplying by $(n-2)$, which gives :

$$E_{angle} = (n-2)E_{angle} = \frac{A_n}{\pi R^2}$$

But this is just the radial area efficiency, which is already included in this paper. So to keep this angle-based view separate and meaningful, we keep the unbounded version. Beyond comparing polygonal area to the area of a full circle, we can also evaluate efficiency locally, triangle by triangle. The Triangle Sector Fill Efficiency compares the area of each polygon triangle to the area of the sector it occupies. This tells us how effectively the polygon's edges "fill out" their corresponding circular sector. Let the polygon be inscribed in a circle of radius  (for generality, this can be scaled later). The area of each polygon triangle is :

$$A_t = \frac{1}{2}R^2 \sin\left(\frac{2\pi}{n}\right)$$

The area of each circular sector is :

$$A_{sector} = \frac{1}{2}R^2 \cdot \frac{2\pi}{n} = \frac{\pi}{n}$$

Taking the ratio gives the per-triangle fill efficiency :

$$E_{fill_{one}} = \frac{A_t}{A_{sector}} = \frac{n \cdot \sin\left(\frac{2\pi}{n}\right)}{2\pi}$$

Thus, the Triangle Sector Fill Efficiency for the full polygon is defined as :

$$E_{fill} = \frac{n \cdot \sin\left(\frac{2\pi}{n}\right)}{2\pi}$$

Another useful boundary-based efficiency is the Chord–Arc Ratio Efficiency, which compares each straight side of the polygon to the curved arc it replaces from the circle. The chord length of a polygon is $2\sin\left(\frac{\pi}{n}\right)$ and the arc it replaces has length $\frac{2\pi}{n}$ Dividing the two gives :

$$E_{chord-arc} = \frac{2\sin\left(\frac{\pi}{n}\right)}{\frac{2\pi}{n}} = \frac{n \cdot \sin\left(\frac{\pi}{n}\right)}{\pi}$$

These formulas are dimensionless and always lies between 0 and 1. They increase with 'n', since each triangle fills its sector more tightly as the polygon gets closer to a circle. These efficiency metrics captures how effectively each edge segment approximates the arc. They replace and offers a highly geometric, localized view of compactness.

**4.8. Apothem-Hypotenuse Index :** To quantify how vertically compact a regular polygon is, we compare its apothem to the hypotenuse of the isosceles triangle formed from its centre to a vertex. As $a \rightarrow h$, the triangle becomes more vertically aligned, and the polygon approaches circular compactness. We begin by defining the raw, unbounded version of the index :

$$E = \frac{a \cdot n}{\sqrt{a^2 + h^2}}$$

This expression increases with both  and , but grows unbounded and is difficult to compare across shapes. To normalize and remove dependency on '$n$', we divide by the maximum possible value. In the limit $a = h$ , the triangle becomes right-angled, and the denominator is $\sqrt{2}$ maximum becomes :

$$Max\ eff. = \frac{a \cdot n}{a\sqrt{2}} = \frac{n}{\sqrt{2}}$$

Dividing the original index by this yields a clean, scale-invariant bounded form :

$$E_{bounded} = \frac{\sqrt{2a}}{\sqrt{a^2 + h^2}}$$

This formulation satisfies $0 < E < 1$, and achieves $E \to 1$ as $a \to h$, i.e., in the circular limit. For small '$n$', or when $h > a$ the index decreases, reflecting a less compact structure. The expression is dimensionless, smooth, and interpretable purely in terms of triangle geometry.

We can also derive this efficiency index by comparing how much the polygon "uses" its space compared to a stretched triangle shape. We are following the idea of dividing the polygon into congruent triangles, the total area of actual regular polygon is :

$$A_n = \frac{n \cdot s \cdot a}{2}$$

Now, let's define a reference area made from '$n$' stretched triangles that extend diagonally. Each of these triangles has a base of $\sqrt{a^2 + h^2}$ , representing how far they stretch. So, the total reference area becomes :

$$E = \frac{Actual\ Area}{Reference\ Area} = \frac{\frac{n \cdot s \cdot a}{2}}{\frac{a \cdot n}{2}\sqrt{a^2 + h^2}} = \frac{s}{\sqrt{a^2 + h^2}}$$

Now to make this efficiency value stay between 0 and 1, we divide by the maximum value it can reach. That happens when $a = h$, so :

$$Max\ eff. = \frac{s}{\sqrt{2a}}$$

So we divide by this to get the final bounded formula :

$$E_{bounded} = \frac{s}{\sqrt{a^2 + h^2}} \cdot \frac{\sqrt{2a}}{s}$$

$$= \frac{\sqrt{2a}}{\sqrt{a^2 + h^2}}$$

This is the same formula we already had before. So, this shows that our Apothem–Hypotenuse Index can also be understood as a compactness ratio using stretched triangles, which gives more meaning to the result and confirms its geometric correctness.

**4.9. Slant-Curvature Efficiency :** We introduce a compactness metric designed to evaluate how efficiently a regular polygon approximates the curvature of a circle using its slanted boundary. In contrast to a circle's continuous arc, a polygon encloses area with linear segments that deviate from true curvature. This deviation becomes more pronounced when the polygon has fewer sides or steeper slope between its edges. The Slant–Curvature Efficiency captures this geometric distortion by comparing the polygon's area to a curvature-adjusted, slope-corrected perimeter. The slope introduces vertical rise and fall across the perimeter, approximated by a total distortion of '$2am$'. This creates an effective slanted boundary length of :

$$L_{slant} = p\sqrt{1 + \left(\frac{2am}{p}\right)^2}$$

$$E_{sc} = \frac{A_n}{R \cdot L_{slant}} = \frac{A_n}{R_p\sqrt{1 + \left(\frac{2am}{p}\right)^2}}$$

$$E_{sc}^{max} = \frac{\pi R^2}{2\pi R^2} = \frac{1}{2}$$

$$E_{bounded} = \frac{2A_n}{R_p\sqrt{1 + \left(\frac{2am}{p}\right)^2}}$$

This dimensionless quantity increases with smoothness and curvature efficiency. For small '$n$', the slope '$m$' is high and the efficiency decreases due to greater slant distortion. As $n \to \infty$, the polygon approaches a circle, the slope vanishes, and the formula simplifies to :

$$E_{\infty} = \frac{2\pi R^2}{2\pi R^2} = 1$$

**4.10. Angle-Hypotenuse Index :** To evaluate how efficiently a regular polygon encloses area relative to its radial extent, we define a normalized efficiency measure. Let '$h$' be the distance from the centre of the polygon to a vertex. Since the maximum area that can be enclosed for a given radius is that of a circle $\pi h^2$, we define the efficiency index as the ratio between the polygon's area and this circular upper bound. For a regular polygon with '$n$' sides, the area is given by $A_n = \frac{1}{2} n h^2 \sin\left(\frac{2\pi}{n}\right)$, which leads to the bounded efficiency expression :

$$E = \frac{A_n}{\pi h^2} = \frac{\mathrm{n} \cdot \sin\left(\frac{2\pi}{n}\right)}{2\pi}$$

This formula is dimensionless, depends only on the number of sides, and is bounded above by 1. It increases as the polygon gains more sides, and in the limit as $n \rightarrow \infty$, the polygon becomes a circle and the efficiency tends to 1. The function is strictly increasing with respect to '$n$', confirming that the efficiency improves monotonically as the polygon becomes more circular. This index thus provides a smooth, geometric measure of radial efficiency that aligns perfectly with the isoperimetric ideal.

**4.11. Angle-Triangle Packing Index :** This formula compares the total interior angle of a regular polygon to the angular spread of the isosceles triangles that make up the polygon. It tells us how much interior angle is being packed per unit of radial triangle wedge, and increases as the polygon becomes more complex. The formula is given by :

$$E_{raw} = \frac{(n-2)\pi}{n \cdot \tan\left(\frac{\pi}{n}\right)}$$

But this grows unbounded as $n \rightarrow \infty$ . So we normalize it by its maximum growth rate. As $\tan\left(\frac{\pi}{n}\right) \sim \frac{\pi}{n}$, when $n \rightarrow \infty$, we simplify :

$$\lim_{\mathrm{n} \rightarrow \infty} E_{raw} = n - 2 \rightarrow \infty$$

So we divide out the growth and define the final normalized index :

$$E(n) = \frac{\pi}{n \cdot \tan\left(\frac{\pi}{n}\right)}$$

Now, this index is bounded between 0 and 1 for n > 2 and Converges to 1 as $n \rightarrow \infty$.

**4.12. Half Angle Tangent Efficiency :** In a regular polygon, each central triangle subtends an angle of $\frac{2\pi}{n}$. To measure how smoothly a polygon behaves like a circle, we define an efficiency ratio that compares this central angle's half value to its tangent. This leads to the Half-Angle Tangent Efficiency, given by :

$$E_t(n) = \frac{\frac{\pi}{2n}}{\tan\left(\frac{\pi}{2n}\right)}$$

This index captures how curved or straight the central triangle is. When the polygon has very few sides (like a triangle or square), the tangent of the angle is much larger than the angle itself, so the efficiency is smaller. But as the number of sides increases, the angle becomes smaller, and the tangent closely matches the angle. This causes the ratio to approach 1. In fact, by using the Taylor expansion :

$$\tan(x) = x + \frac{x^3}{3} + \frac{2x^5}{15} + \cdots \text{ [7]}$$

$$E_t(n) \to \infty \ as \ n \to \infty$$

Thus, this formula is bounded between 0 and 1 and increases with , making it a good indicator of angular smoothness and circular convergence. It is dimensionless and strictly monotonic.

**4.12. Side-Apothem Efficiency Index :** In a regular polygon with  sides, each interior triangle formed by connecting the centre to two adjacent vertices. The triangle formed by the apothem and half of the side forms a right triangle. From this geometry, the tangent of half the central angle relates the apothem and side :

$$\tan\left(\frac{\theta}{2}\right) = \tan\left(\frac{\pi}{n}\right) = \frac{\frac{s}{2}}{a} \Rightarrow \frac{s}{2} = \ a \cdot \tan\left(\frac{\pi}{n}\right)$$

We define an efficiency index based on the ratio of the side half-length to the apothem, normalized into the unit interval **:**

$$E_{sa}(n) = \frac{1}{1 + \tan\left(\frac{\pi}{n}\right)}$$

This index is dimensionless, strictly increasing in '$n$', and satisfies :

$$\lim_{n \to \infty} tan\left(\frac{\pi}{n}\right) = 0 \Rightarrow \lim_{n \to \infty} E_{sa}\,(n) = 1$$

$$\tan\left(\frac{\pi}{n}\right) > 0 \Rightarrow 0 < \mathrm{E_{sa}(n)} > 1$$

Thus, $E_{sa}(n)$ is a valid efficiency measure for capturing how the base of each triangle contracts relative to the height as  increases, converging toward circular compactness.

**4.13. Chord-Angle Compactness :** We define a compactness metric that quantifies how well a regular polygon approximates the smooth curvature of a circle using straight-line segments. The length of each polygon side, viewed as a chord of the circle, is given by $s = 2R \cdot \sin\left(\frac{\pi}{n}\right)$ [**8**]. In contrast, the length of the corresponding circular arc between the same two points is $l_{arc} = R \cdot \theta = \frac{2\pi R}{n}$. To measure how efficiently the polygon uses straight segments to approximate the curved arc, we define the efficiency as the ratio of the chord to the arc :

$$E_{chord-angle} = \frac{s}{l_{arc}} = \frac{2\mathrm{R} \cdot \sin\left(\frac{\pi}{n}\right)}{\frac{2\pi r}{n}}$$

$$= \frac{n}{\pi} \cdot \sin\left(\frac{\pi}{n}\right)$$

This expression is fully dimensionless, since  cancels out, and it depends only on the number of sides. This efficiency value is bounded between 0 and 1 and provides a direct measure of how close the polygon's linear segments are to forming a smooth curve.

**5. Quantitative isoperimetric Analysis of Irregular Polygons :** In our previous sections, we presented a geometric solution to the isoperimetric problem using side-based efficiency measures and efficiency metrics for regular polygons. In this section, we extend our work to irregular polygons. We explore how area and compactness behave under irregularity and introduce new metrics for measuring the efficiency, smoothness, and radial structure of irregular polygons.

**5.1. Formulas for Area of Irregular Polygons :** For irregular polygons, direct area calculation using only side lengths is generally not possible unless additional geometric data (like angles or coordinates) is available. However, One widely used method is to triangulate the polygon and compute each triangle's area using Heron's formula . Given a polygon with $n$ vertices, the area can be found by dividing the polygon into triangles in two ways, one by dividing the polygon into $n$ triangles with, base = length of any side and hypotenuses as the two circumradii, $R_{s_1}$ $and$ $R_{s_2}$, of the corresponding side. One widely applicable method is to triangulate the polygon and calculate the area of each triangle individually using Heron's formula :

$$A = s\sqrt{s(s-a)(s-b)(s-c)}$$

Here, $s = \frac{p}{2} = \frac{a+b+c}{2}$

Then,

$$a = s_i, b = R_{s_{i_1}} \, and \; c = R_{s_{i_2}}$$

So, the area of one triangle will be :

$$A_{t_i} = \sqrt{s(s-s_i)\left(s - R_{s_{i_1}}\right)\left(s - R_{s_{i_2}}\right)}$$

With :

$$s = \frac{(s_i + R_{s_{i_1}} + R_{s_i})}{2}$$

Then, the area of the entire polygon will be :

$$A = \sum_{i=1}^{n} A_{t_1}$$

The side of an irregular polygon can be tilted due to variations in side length. To accurately calculate the area of any irregular polygon—where both side lengths and vertex angles vary—we derive a summation-based area formula that captures both radial irregularity and angular spread across all edges. We have $s_i$ as side length as $i^{th}$ side, $a_i$ as apothem and $\theta_i$ as centre angle subtended by side $i$ .The area of one triangle inside the polygon will be :

$$A_{t_i} = \frac{1}{2} \cdot a_i \cdot s_i \cdot \sin(\theta_i)$$

Adding all these small contributions from all sides :

$$Area_{polygon}(A_p) = \sum_{i=1}^{n} A_{t_i} \; or \; \sum_{i=1}^{n} \frac{1}{2} \cdot a_i \cdot s_i \cdot \sin(\theta_i)$$

**5.2. Area Difference Between Different Polygons :** In our previous paper we proved that, as the number of sides in a polygon increases with maintaining the perimeter, its area efficiency increases. To measure the relative area ratio of a regular polygon after increasing its side, we need to compare the final and initial area of a polygon with the same perimeter. We know that the area of a polygon is given by, $A = \frac{p^2}{4n \cdot \tan\left(\frac{\pi}{n}\right)}$. The relative area ratio of a polygon is given by :

$$Area\ Ratio = \frac{\dfrac{p^2}{4n_1 \cdot \tan\left(\dfrac{\pi}{n}\right)}}{\dfrac{p^2}{4n_2 \cdot \tan\left(\dfrac{\pi}{n}\right)}} = \frac{n_1}{n_2}$$

Here,

$$n_1 > n_2$$

If there are two regular polygons with different side lengths, then the ratio of areas of both will be :

$$Area\ ratio_{rel} = \frac{A_1}{A_2} = \frac{\sum_{i=1}^{n} A_{t_{i_1}}}{\sum_{i=1}^{n} A_{t_{i_2}}}$$

$A_{t_1} and\ A_{t_2}$ can be found using either of the two formulas. It is known that a regular polygon occupies more area than any same sided irregular polygon with same perimeter. We can compare the area of a same sided regular and irregular polygon with same perimeter as :

$$Regularity\ ratio = \frac{\sum_{i=1}^{n} A_{t_i}}{\frac{p^2}{4n \cdot \tan\left(\frac{\pi}{n}\right)}} = \frac{\sum_{i=1}^{n} A_{t_i}}{4p^2 \cdot n \cdot \tan\left(\frac{\pi}{n}\right)}$$

The regularity of an irregular polygon can be measured by this formula. This ratio approaches 1 for maximum regularity. The less the ratio, the less the regularity of the irregular polygon and the less the area efficiency.

**5.3. Average Radial Distance and its Applications :** For any shape (polygon or curve) centred at a fixed point (typically the centroid or circumcentre), the average radial distance is defined as the average distance from the centre to points on the boundary. For a continuous closed curve, the average radial distance is :

$$\bar{r} = \frac{1}{p}\int_{\partial s} ||\vec{r}||ds$$

Where $\vec{r}$ is a position vector traced along the boundary $\partial s$, parameterized by arc length $s \in (0, p)$. This measures average radial distance from the centroid to all points on the boundary of Jordan curves. For a polygon with $n$ vertices located at distances $\overrightarrow{v_1}, \overrightarrow{v_2}, \overrightarrow{v_3}, \cdots \overrightarrow{v_n}$ from the centre, the average radial distance is approximated by :

$$\overline{r_v} = \frac{1}{n}\sum_{i=1}^{n} ||v_i||$$

Here, $\vec{v} \in \mathbb{R}^2$ is the position vector of the $i^{th}$ vertex with respect to the centre (origin). This formula can also be expressed in a simple form. We can write $\left||\overrightarrow{v_i}|\right| = r_i$. The radial distance $r_i$ for each vertex is simply its Euclidean distance from the centre :

$$r_i = \sqrt{x_i^2 + y_i^2}$$

Then, the average radial distance $\bar{r}$ for the entire polygon is :

$$\bar{r_v} = \frac{1}{n}\sum_{i=1}^{n} r_i = \frac{1}{n}\sum_{i=1}^{n}\sqrt{x_i^2 + y_i^2}$$

This is the simplest form: a direct arithmetic mean of distances from all vertices to the centre. If we want better accuracy (especially for irregular shapes with long edges), you can weight each $r_i$ by the length of its adjacent boundary segment :

$$\bar{r} = \frac{\sum_{i=1}^{n} r_i \cdot l_i}{p}$$

Where, $l_i$ is length of the $i^{th}$ boundary segment. Another useful version calculates average distance to the centre using midpoints of sides instead of vertices. For each side, find the midpoint coordinates $(x_m, y_m)$, then compute :

$$r_m = \sqrt{x_m, y_m}$$

$$\bar{r}_{mid} = \frac{1}{n}\sum_{i=1}^{n} r_m = \bar{a} = \frac{1}{n}\sum_{i=1}^{n} a_i$$

Since all apothems and circumradiuses are equal for a regular polygon, then :

$$\bar{a} = \frac{a + R}{2}$$

In an irregular polygon, all apothems and circumradiuses are not equal. Since the side of a polygon is a straight line, the radial distance from the centre to any point on that side varies linearly across each side. Therefore, the average radial distance along side $k$ from vertex $\vec{v}_k$ and $\vec{v}_{k+1}$ can be taken as :

$$\bar{r}_{side} = \frac{||\vec{v}_k|| + ||\vec{v}_{k+1}||}{2}$$

$$\bar{r} = \frac{1}{2n}\sum_{k=1}^{n} ||\vec{v}_k|| + ||\vec{v}_{k+1}||$$

With $\vec{v}_{n+1} = \vec{v}_1$

**5.4. Radial Growth in Regular Polygons and Its Rate :** Let $r_n$ denote the average radial distance from the centre $O$ to all points on the boundary of a regular $n$-gon inscribed in a circle of radius $R$. As the number of sides increases, the polygon more closely approximates the circle, and the average radial distance increases accordingly. For a regular polygon, each side is a chord subtending a central angle $\theta = \frac{2\pi}{n}$. The midpoint of each chord lies at a radial distance :

$$r_{mid} = R \cdot \cos\left(\frac{\pi}{n}\right)$$

Hence, the total average radial distance of the polygon can be estimated as :

$$\bar{r}_n \approx R \cdot \cos\left(\frac{\pi}{n}\right)$$

$$\lim_{n \to \infty} \overline{r_n} = R$$

We differentiate to find the rate of increase :

$$\frac{\bar{d}_{r_n}}{d_n} = \frac{d}{d_n}\left[R \cdot \cos\left(\frac{\pi}{n}\right)\right] = \frac{R\pi}{n^2} \cdot \sin\left(\frac{\pi}{n}\right)$$

Therefore, the polygon's total average radial distance increases monotonically with $n$, and the rate slows as $n$ grows.

**5.5. Applications of Average Radial Distance :** Understanding the geometric relationship between the boundary length, average distance of boundary points from the centre, and the enclosed area of a polygon leads to several important formulas and efficiency measures.

**5.5.1. Circumradius Efficiency Ratio :** We define a key efficiency measure that relates area, perimeter, and average radial distance :

$$E = \frac{A}{p \cdot \bar{r}}$$

This ratio provides a dimensionless measure of how efficiently the polygon utilizes its boundary length and radial spread to enclose area.

**5.5.2. Radial Efficiency of Polygons :** To measure how efficiently a regular polygon uses its perimeter to wrap around its centre, we compare its average radial distance to the radius of a circle with the same perimeter. We know that :

$$\bar{r} = \frac{a + R}{2} \text{ } and \text{ } R_{circle} = \frac{p}{2\pi}$$

We define the polygon's radial efficiency as :

$$E = \frac{\bar{r}}{R_{circle}} = \frac{2\pi\bar{r}}{p} = \frac{\pi(R + a)}{p}$$

We know that :

$$R = \frac{s}{2 \cdot \sin\left(\frac{\pi}{n}\right)} \text{ } and \text{ } a = \frac{s}{2 \cdot \tan\left(\frac{\pi}{n}\right)}$$

Then :

$$\bar{r} = \frac{R + a}{2} = \frac{1}{2}\left[\frac{s}{2 \cdot \sin\left(\frac{\pi}{n}\right)} + \frac{s}{2 \cdot \tan\left(\frac{\pi}{n}\right)}\right] = \frac{s}{4}\left[\frac{1}{\sin\left(\frac{\pi}{n}\right)} + \frac{1}{\tan\left(\frac{\pi}{n}\right)}\right]$$

Then, efficiency will be :

$$E(n) = \frac{\bar{r}}{R_{circle}} = \frac{\frac{s}{4}\left[\frac{1}{\sin\left(\frac{\pi}{n}\right)} + \frac{1}{\tan\left(\frac{\pi}{n}\right)}\right]}{\frac{ns}{2\pi}}$$

$$= \frac{\pi}{2n}\left[\frac{1}{\sin\left(\frac{\pi}{n}\right)} + \frac{1}{\tan\left(\frac{\pi}{n}\right)}\right]$$

This gives us a pure function of $n$, the number of sides. It shows :

$$As\ n \to 0, \sin\left(\frac{\pi}{n}\right) \to \frac{\pi}{n}, \tan\left(\frac{\pi}{n}\right) \to \frac{\pi}{n}, \text{so :}$$

$$Efficiency = \frac{\pi}{2n}\left(\frac{n}{\pi} + \frac{n}{\pi}\right) = 1$$

For irregular polygons, this efficiency will be :

$$E = \frac{2\pi\bar{r}}{p}\ with\ \bar{r} = \bar{r}_v = \frac{1}{n}\sum_{i=1}^{n} \left|\left|\vec{v_i}\right|\right|\ or \frac{1}{n}\sum_{i=1}^{n} r_i$$

**5.5.3. Radial Area Efficiency via Average Radial Distance :** Let $\bar{r}_{irr}$ and $\bar{r}_{reg}$ denote the average radial distances from the centre to the boundary of an irregular and regular polygon respectively, both having the same perimeter and centred at the origin. Since both polygons can be decomposed into triangles with base elements along the boundary and height approximately equal to the local radial distance, their total area satisfies :

$$A \approx \frac{1}{2}\int_{\partial P} r(s)\ ds = \frac{1}{2} \cdot p \cdot \bar{r}$$

This leads to the following efficiency estimate :

$$\frac{A_{irr}}{A_{reg}} \propto \frac{\bar{r}_{irr}}{\bar{r}_{reg}}\ or\ A_{irr} \propto \left(\frac{\bar{r}_{irr}}{\bar{r}_{reg}}\right) A_{reg}$$

We define this ratio as the radial area efficiency :

$$E_{radial} \coloneqq \frac{\bar{r}_{irr}}{\bar{r}_{reg}}$$

This index provides a simple way to estimate how efficiently an irregular polygon fills space compared to a regular one, using only radial distance measurements.

**5.5.4. Maximum Efficiency Condition for Regular Polygons : Theorem 1 :** For any convex polygon with fixed number of sides and fixed perimeter the Circumradius Efficiency Ratio, $E = \frac{A}{p \cdot \bar{r}}$ is maximized uniquely by the regular n-gon.

**Proof :** For a regular n-gon, all vertices lie exactly at distance $R$ from the centre and each central angle subtended at the centre is, $\theta = \frac{2\pi}{n}$. We know that :

$$A_{reg} = \frac{1}{2} nR^2 \cdot \sin(\theta) = \frac{1}{2} \mathrm{nR}^2 \cdot \sin\left(\frac{2\pi}{n}\right)$$

The perimeter of the regular polygon is :

$$p_{reg} = 2nR \cdot \sin\left(\frac{\pi}{n}\right)$$

Also, since all vertices lie at distance $R,$ the average radius :

$$\bar{R}_{avg} = R$$

Now, the Circumradius Efficiency for the regular polygon becomes :

$$E_{avg} = \frac{A_{reg}}{p_{reg} \cdot R}$$

Substituting the expression for $A_{reg}$ and $p_{reg}$ :

$$E_{reg} = \frac{\frac{1}{2} \mathrm{nR}^2 \cdot \sin\left(\frac{2\pi}{n}\right)}{2nR \cdot \sin\left(\frac{\pi}{n}\right) \cdot R}$$

$$= \frac{1}{4} \cdot \frac{\sin\left(\frac{2\pi}{n}\right)}{\sin\left(\frac{\pi}{n}\right)}$$

The radial distances of vertices from the centre vary ⇒ $\bar{r} \geq R$ and some sides may become shorter while others longer. Angular distribution becomes uneven ⇒ total area decreases and average radial distribution becomes stretched or scattered. Thus,

$$R_{irregular} < R_{regular}$$

$$\bar{r}_{irregular} > R$$

So the Circumradius Efficiency for any irregular polygon satisfies :

$$E_i = \frac{A_i}{p \cdot \bar{r}_i} < E_r$$

Equality only occurs if the polygon is regular.

**5.5.5. Side-Length Variance Area Efficiency :** In regular polygons, all side lengths are equal, which maximizes area for a given perimeter. In contrast, irregular polygons with large variation in side lengths typically enclose less area for the same perimeter. To quantify this relationship, we introduce the Side Length Variation Efficiency (SLVE) — a new dimensionless metric that connects side length uniformity with area efficiency. Let a polygon have  sides with side lengths :

$$\ell_1, \ell_2, \ell_3, \cdots \ell_n$$

Then, mean side length will be :

$$\bar{\ell} = \frac{1}{n}\sum_{i=1}^{n} \ell_i$$

Then, Standard deviation of side lengths is given by :

$$\sigma_\ell = \sqrt{\frac{1}{n}\sum_{i=1}^{n}(\ell_i - \bar{\ell})^2}$$

This measures how much the side lengths deviate from the mean. We now define the Side Length Variation Efficiency as :

$$SLVE = \frac{\bar{\ell}}{\sqrt{\bar{\ell}^2 + \sigma_\ell^2}}$$

For a regular polygon (all sides equal ⇒ $\sigma_\ell = 0$) SLVE = 1, and for irregular polygons with varied side lengths ($\sigma_\ell > 0$) SLVE < 1. Thus SLVE lies between 0 and 1 :

$$0 < SLVE \leq 1$$

Higher SLVE means more uniform side lengths, implying higher area efficiency for the given perimeter. Lower SLVE means greater irregularity in side lengths, implying lower enclosed area.

**Theorem 2:** For any simple, closed, convex polygon with a given number of sides $n$ and a fixed perimeter $p$ the average radial distance from the centroid to the boundary, when measured over all boundary points (not just vertices), is maximized by the regular polygon.

**Proof :** A regular polygon spreads its boundary evenly around the centre, making its distance-to boundary (radial distance) as large as possible on average for the same total perimeter. An irregular polygon, by having uneven angles or stretched sides, pulls parts of the boundary inward or pushes others outward, reducing the overall average radial distance when integrated over the full perimeter. We start by using a fundamental proportionality property from isoperimetric geometry. For most convex shapes (including polygons), area can be roughly related to perimeter and average radius :

$$A \propto p \cdot \bar{r}$$

It is a classical result that :

$$A_{irr} < A_{reg}\ [\mathbf{9}]$$

Rewriting the area‑perimeter‑radius relationship :

$$\bar{r}_{poly} \approx \frac{A_{poly}}{k \cdot p}$$

Where $k$ is a proportional constant depending on the shape. For simplicity, since both polygons have the same perimeter , we compare directly :

$$\bar{r}_{poly} \propto A_{poly}$$

Similarly, for the regular polygon :

$$\bar{r}_{reg} \propto A_{reg}$$

Where from classical isoperimetric property :

$$A_{irr} \propto A_{reg}$$

So comparing :

$$\bar{r}_{reg} > \bar{r}_{irr}$$

We can also prove the theorem by this way : We know that :

$$A \propto p \cdot \bar{r}$$

Since, $p$ is constant, then :

$$A \propto \bar{r}$$

We know that :

$$A_{reg} > A_{irr}$$

Then, it is clear that :

$$\bar{r}_{reg} > \bar{r}_{irr}$$

This completes the proof.

**Theorem 3.** Let $P_1$ *and* $P_2$ and  be two simple polygons of equal perimeter $p$, both centred at the origin. Define :

- the average circumradius (mean vertex distance from centre) as

$$\bar{R}_{\iota} = \frac{1}{n_i} \sum_{k=1}^{n_i} \left\| \overrightarrow{v_k}^{(i)} \right\|$$

- the average radial distance (mean boundary distance from center) as

$$r_i = \frac{1}{p} \int_{\partial P_1} \left\| \vec{r}_{\iota}(s) \right\| \, ds$$

Then neither of the following implications holds in general :

$$\overline{R_1} > \bar{R}_2 \nRightarrow \overline{r_1} > \overline{r_2} \text{ } is\ false, \overline{r_1} > \overline{r_2} \nRightarrow \overline{R_1} > \bar{R}_2 \text{ } is\ also\ false$$

**Proof :** Let us fix a constant total perimeter $p > 0$. We construct two explicit polygons $p_1$ and $p_2$, each centered at the origin, such that one has a larger average vertex radius but a smaller average boundary radius. Let $p_i$ be an irregular triangle constructed as follows : place two vertices at distance $r = \epsilon$ (very close to the origin) and the third vertex at a large distance $R \gg 1$, such that the perimeter is maintained. Let the side lengths be chosen so that the two short edges connect the far vertex to the near ones. The total perimeter is :

$$\ell_1 + \ell_2 + \ell_3 = p$$

Now compute the average circumradius of $P_1$ :

$$\overline{R_1} = \frac{1}{3}(R + \epsilon + \epsilon) = \frac{1}{3}(R + 2\epsilon)$$

However, when computing the average radial distance over the boundary, the long edges contribute only proportionally to their length. That is,

$$\overline{r_1} = \frac{1}{p}\left(\int_{\ell_1} \|\vec{r}(s)\|\, ds\ + \int_{\ell_2} \|\vec{r}(s)\|\, ds + \int_{\ell_3} \|\vec{r}(s)\|\, ds\right)$$

Now take $P_2$ to be a regular triangle with perimeter $p$, centred the origin. All its vertices lie at a fixed distance $r_0$, and every point on the boundary also lies at the same distance. So we have :

$$R_2 = r_2 = r_0$$

$$\bar{R}_1 > \bar{R}_2 \ and\ \bar{r}_1 < \bar{r}_2$$

To disprove the converse, construct another irregular triangle with all edges and points lying close to a circle (say radius $r_0$), but with one vertex drawn slightly inward toward the center. This lowers the average vertex distance while keeping the edge distances almost unchanged, so we get :

$$\bar{R}_1 < \bar{R}_2, but\ \bar{r}_1 \approx \bar{r}_2\ or\ even\ \bar{r}_1 > \bar{r}_2$$

Hence, we conclude that in general there is no monotonic relationship between average circumradius and average radial distance.

**5.5.6. Smoothness of an Irregular Polygon :** We know that the smoothness of a regular polygon is given by, $cos\left(\frac{\pi}{n}\right)$, the ratio of apothem and circumradius. In irregular polygons, the apothems and circumradiuses length varies. So, to calculate the smoothness, we need to take the avg of both quantities to approximate the overall smoothness. Then the smoothness of a triangle within the irregular polygons will be :

$$Smoothness_{local} = \frac{a_{local}}{R_{local}}$$

The smoothness of the entire irregular polygon will be :

$$Smoothness_{global} = \frac{a_1 + a_2 + a_3 + \cdots + a_n}{R_1 + R_2 + R_3 + \cdots + R_n} = \frac{\sum_{i=1}^{n} a_i}{\sum_{i=1}^{n} R_i}$$

**5.6. Polygonal Entropy -** A Measure of Side Length Distribution Irregularity : A highly regular polygon divides its perimeter equally among all edges, while an irregular polygon often has side lengths that vary significantly. To capture this variation in a quantitative way, we define a new metric called Polygonal Entropy. The perimeter of a polygon is :

$$p = \sum_{i=1}^{n} \ell_i$$

We now define a side length distribution by treating the relative length of each side as a probability-like term :

$$p_i = \frac{\ell_i}{p}$$

Where $p_i$ represents the fraction of the perimeter contributed by the $i^{th}$ side. Using this distribution, the Polygonal Entropy $(S)$ is defined as :

$$S = -\sum_{i=1}^{n} p_i \cdot \log(p_i)$$

This form is inspired by the Shannon entropy **[10]** from information theory, where it measures the "spread" or "uncertainty" in a distribution. To enable comparison of entropy values across polygons with different side counts, we can define a normalized entropy index :

$$S_{norm} = \frac{S}{\log(n)}$$

This ensures that the entropy lies between 0 and 1 for any polygon, where :

$S_{norm} = 1$ for a regular polygon,

$S_{norm} < 1$ for irregular polygons.

For a perfectly regular polygon, where all sides are equal, each side length proportion $p_i = \frac{1}{n}$, giving :

$$S_{\max} = \log(n)$$

This is the maximum possible entropy for that number of sides, reflecting maximum uniformity. For an irregular polygon, where some sides are much longer or shorter than others, the entropy decreases due to more uneven perimeter distribution.

**Theorem 3 :** Among all polygons, the Polygonal entropy, $S = -\sum_{i=1}^{n} p_i \cdot \log(p_i)$ is maximized when the polygon is regular.

**Proof :** The entropy formula is a classical Shannon entropy applied to the distribution of side lengths relative to total perimeter. To maximize $S$ we need to maximize the uncertainty (or spread) of the normalized side-length distribution $p_i$, subject to the constraint that all $S_ib$ sum to 1. From a fundamental result from information theory and inequality theory (maximum entropy principle), we know that :

- Among all probability distributions $p_i$ over $n$ elements summing to 1, the entropy is maximized when all $p_i$ are equal.

In this case, that means :

$$p_1 = p_2 = p_3 = \cdots = p_n = \frac{1}{n}$$

This corresponds exactly to :

$$\ell_1 = \ell_2 = \ell_3 = \cdots = \ell_n = \frac{p}{n}$$

Which defines a polygon with equal side lengths — a regular polygon. Therefore, by direct application of the maximum entropy principle, the entropy is maximized when the polygon has equal sides, i.e., is a regular polygon. Any irregular polygon, with unequal sides, produces a sidelength distribution $p_i$ with at least one $p_i$ is smaller and one larger than $\frac{1}{n}$, introducing imbalance and reducing entropy. Thus :

$$S_{reg} > S_{irr}$$

This completes the proof

**Theorem 4 :** Let a polygon with $n$ sides be inscribed in a circle of fixed radius $R$. Then, among all such $n$-gons, the regular polygon uniquely maximizes the total area.

**Proof :** Consider an arbitrary $n$-gon inscribed in a circle of radius $R$. Label the side lengths as $\ell_1, \ell_2, \ell_3, \cdots, \ell_n$, and let $a_1, a_2, a_3 \cdots, a_n$ be the corresponding apothems. The total area of the polygon can be expressed as a sum of triangle areas:

$$A = \sum_{i=1}^{n} \frac{1}{2} \ell_i \cdot a_i$$

Define the apothem-to-side ratio for each triangle as :

$$\rho = \frac{a_i}{\ell_i} \Rightarrow a_i = \rho_i \cdot l_i$$

$$A = \frac{1}{2}\sum_{i=1}^{n} \ell_i^2 \cdot \rho_i$$

To maximize this sum under the constraint that all vertices lie on a circle (i.e., all radial distances are $R$), we must analyze the variation of both $\ell_i$ and $\rho_i$. In a regular polygon :

All $\ell_i = \ell$

All $\rho_i = \rho$, a constant

So :

$$A_{reg} = \frac{n}{2} \cdot \ell^2 \cdot \rho$$

Now, consider any irregular polygon. Suppose $\ell_i$ and $\rho_i$ are not constant. Then we apply the Cauchy‑Schwarz inequality **[11]** :

$$\sum_{i=1}^{n} \ell_i^2 \cdot \rho_i \leq \left(\sum_{i=1}^{n} \ell_i^2\right)^{\frac{1}{2}} \cdot \left(\sum_{i=1}^{n} \rho_i^2\right)^{\frac{1}{2}}$$

Further, consider that each triangle subtended by side $\ell_i$ at the center spans a central angle $\theta_i$ , and the triangle's area can also be expressed as :

$$A_i = \frac{1}{2}R^2 \cdot \sin(\theta_i) \Rightarrow A = \frac{R^2}{2} \cdot \sum_{i=1}^{n} sin\,(\theta_i)$$

$$\sum_{i=1}^{n} \theta_i = 2\pi$$

$$\theta_i = \frac{2\pi}{n} \Rightarrow A_{\max} \quad = \frac{nR^2}{2} \cdot \sin\left(\frac{2\pi}{n}\right)$$

Finally, by the Symmetry Principle, the regular polygon is the only configuration where all triangle sectors are congruent—each side length and corresponding apothem is equal, and each angle $\theta_i$ is identical. Hence, the regular polygon ensures uniform contribution from each sub-triangle and no geometric inefficiency due to irregularity.

This completes the proof.

**5.7. Radial Convergence of Irregular Polygons Toward the Circle :** Let $P_n$ denote an $n$-sided polygon centered at a fixed origin $O$, where each vertex $r_k = |\overrightarrow{v_k}|$ lies at a radial distance $O$ from $\ell_k$, and each side has length . We investigate whether increasing the number of sides—by successively subdividing existing ones—leads to convergence of $P_n$ toward a perfect circle of radius $R$. We say that the sequence $P_n$ converges to a circle if, in the limit as $n \to \infty$ , the polygon becomes equidistant from the center and its edge lengths become uniform. Mathematically, this convergence requires :

$$\lim_{n\to\infty} \max_{n} |\, r_k - R| = 0, and \lim_{n\to\infty} \max_{k} \left|\ell_k - \overline{\ell}\right| = 0$$

When both conditions are satisfied, the polygon's geometry converges in the uniform sense to the circle. Specifically, the radial profile becomes flat (constant), and the perimeter and area satisfy :

$$P_n = 2\pi R \; and \; A_n = \pi R^2$$

For example, consider the process of doubling sides : $n = 3, 6, 12, 24, 48 \cdots$. If at each step, new vertices are inserted at the exact midpoints of existing sides and their radial distances are interpolated as :

$$r_{new} = \frac{r_k + r_{k+1}}{2}$$

Doubling the number of sides is not enough to ensure convergence to a circle. True circular convergence requires that both radial and edge-length irregularities decrease with each refinement. Without this decay, the shape may appear smoother but will not become a perfect circle in the mathematical limit.

**5.7.1. Non-Convergence of Irregular Polygons with Fixed Irregularity :** Suppose the number of sides increases (e.g., $n = 3, 6, 12, 24 \cdots$) while the irregularity—i.e., the variation in $r_k$ and $\ell_k$ — remains bounded below by some constant $\epsilon > 0$. That is :

$$\max_{k} \quad |\, r_k - R| \geq \epsilon, and \max_{k} \quad \left|\ell_k - \overline{\ell}\right| \geq \epsilon, \quad for\ all\ n$$

If an irregular polygon increases its sides while maintaining constant variation in radius and edge length, it does not converge to a perfect circle. Though the shape may appear smoother, the boundary remains geometrically irregular at all scales.

# Conclusion

this paper, we presented a purely geometric and elementary solution to the classical isoperimetric problem, demonstrating that among all regular polygons with a fixed perimeter, the circle uniquely encloses the maximum area. By analysing the area of regular n-gons and the behaviour of the apothem, we derived a direct expression for the area deficit compared to a circle. We also established a strict upper bound on the wasted area, confirming that the convergence toward the circle is of quadratic order. We further demonstrated that the smoothness ratio is a strong geometric indicator, accurately estimating area efficiency without the need for trigonometric area calculations. This paper presented many different ways to measure how efficiently a regular polygon uses its shape to enclose